\documentclass{ifacconf}
\pdfoutput=1
\usepackage{natbib,amsmath,amssymb,amsfonts}
\usepackage{graphicx}      
\usepackage{microtype}
\newcommand{\R}{\mathbb{R}}

\newcommand{\cond}{\operatorname{cond}}
\newcommand{\cper}{C_\mathrm{per}}
\newcommand{\avg}{\mathrm{avg}}
\renewcommand{\d}{\,\mathrm{d}}
\begin{document}
\begin{frontmatter}

  \title{Using feedback control and Newton iterations to track
    dynamically unstable phenomena in experiments}


\author[First]{Jan Sieber} 
\author[Second]{Bernd Krauskopf} 

\address[First]{Department of Mathematics, University of Portsmouth,
  Lion gate Building, Lion Terrace, Portsmouth, PO1 3HF, United
  Kingdom.}  
\address[Second]{Department of Engineering Mathematics,
  University of Bristol, University Walk, Bristol, BS8 1TR, United Kingdom.}

\begin{abstract}                
  If one wants to explore the properties of a dynamical system
  systematically one has to be able to track equilibria and periodic
  orbits regardless of their stability. If the dynamical system is a
  controllable experiment then one approach is a combination of
  classical feedback control and Newton iterations. Mechanical
  experiments on a parametrically excited pendulum have recently shown
  the practical feasibility of a simplified version of this algorithm:
  a combination of time-delayed feedback control (as proposed by
  Pyragas) and a Newton iteration on a low-dimensional system of
  equations. We show that both parts of the algorithm are uniformly
  stable near the saddle-node bifurcation: the experiment with
  time-delayed feedback control has uniformly stable periodic orbits,
  and the two-dimensional nonlinear system which has to be solved to
  make the control non-invasive has a well-conditioned Jacobian.
\end{abstract}


\end{frontmatter}

\section{Introduction}
\label{sec:intro}
One way to explore a nonlinear dynamical system in a systematical
fashion is \emph{bifurcation analysis} by continuation: one starts in
a parameter region where one knows that a simple attractor exists
(say, a stable periodic orbit) and then varies a system parameter $p$,
checking at which parameter values the periodic orbit loses its
stability or ``disappears''. At these special parameter values the
periodic orbit undergoes a \emph{bifurcation}, and other invariant
objects (equilibria, periodic orbits, tori) may branch off. Thus, by
systematically tracking equilibria, periodic orbits and their
bifurcations in the parameter space one can (to a good extent)
classify the qualitative behaviour of the dynamical system.

In this paper we analyse the stability of the algorithm for
pseudo-arclength continuation in experiments introduced in
[\cite{SGNWK08}] that considered periodic rotations of a vertically
excited pendulum near a saddle-node bifurcation.  For this example
experiment we demonstrate (using simulations) that all parts of the
algorithm converge uniformly. In section~\ref{sec:parc} we recap how
one embeds a Newton iteration into pseudo-arclength continuation to
study how periodic orbits depend on system parameters. We also briefly
explain how this continuation is implemented for the excited pendulum
using time-delayed feedback control (TDFC; \cite{P92}) for the
periodic part of the problem and a Newton iteration for only two
scalar variables.  We show that, with this approach, the system with
TDFC is uniformly stable near the saddle-node bifurcation (which is
not the case for the classical TDFC), and that the Jacobian used in
the Newton iteration is uniformly well-conditioned.

\section{Background on related methods}
\label{sec:back}
If a dynamical system is given in the form of a low-dimensional
ordinary differential equation (or discrete map) one can apply
specialized algorithms based on Newton iterations embedded into
pseudo-arclength continuation (see Section~\ref{sec:parc} for an
explanation and an example), which are available as software packages,
for example, \textsc{Auto} or \textsc{Matcont} (see textbooks
[\cite{D07,K04}] for a detailed introduction). These algorithms have
been successfully extended to problems where the models are delay
differential equations (\textsc{Dde-Biftool}: \cite{ELS01};
\textsc{Pddecont}: \cite{SSH04}), dissipative partial differential
equations (\textsc{Loca}: \cite{LOCA02}; \cite{LRSC98}), or
high-dimensional systems (such as stochastic Monte Carlo simulations)
with `essentially low-dimensional' dynamics (\cite{KGH04}). One
advantage of algorithms based on Newton iterations is that they work
independent of the dynamical stability of the state they track. Hence,
they are also able to track unstable periodic motions or equilibria.

In contrast to the situation where one investigates a model,
bifurcation analysis in experiments is typically done by parameter
studies: one gradually varies a system parameter and observes the
transients. When one observes a slow-down of the transients or a
sudden `jump' of the output it is likely that one has encountered a
bifurcation (and its type can sometimes be deduced from the transient
behavior). Using an electronic implementation of a Duffing oscillator,
\cite{LP02} show how one can automate this approach to trace out
bifurcations. \cite{ASFKRK99} for an eloctrochemical system and
\cite{DM03} for a Colpitts oscillator have demonstrated the use of
Newton iterations and pseudo-arclength continuation in experiments
(also detecting or continuing bifurcations). Both studies ran a system
identification procedure in parallel to the experiment, determined the
steady states (fixed points or periodic orbits) of the identified
model, and used feedback control to drive the experiment toward the
identified steady state. The accuracy of the results using this
approach is limited not only by the measurement accuracy and the
tolerances set in the Newton iteration but also by the accuracy of the
system identification. This is a severe handicap because system
identification is an inverse and, thus, ill-posed
problem. \cite{ASFKRK99} demonstrated their approach for simulations
of the chemical system only, actual experiments are still outstanding.

An alternative approach to continuation of steady states in
experiments is via the use of feedback control mechanisms that are
automatically \emph{non-invasive}: washout filters [\cite{AWC94}] and
time-delayed feedback [TDFC; \cite{P92,SSG94,KBPOMBRE01}]. Whereas
classical feedback control compares the output $\phi$ of a dynamical
system to a given reference signal $\tilde\phi$ and feeds a (typically
linear) combination of the difference $\phi-\tilde\phi$ back into the
experiment, washout-filtered feedback control and TDFC do not require
a given reference signal. Washout filtered feedback picks $\tilde\phi$
as the solution of
\begin{displaymath}
  \Dot{\tilde\phi}=R\cdot(\phi-\tilde\phi)  
\end{displaymath}
where $R$ is a stable matrix. TDFC picks $\tilde\phi$ using the
recursion
\begin{equation}
  \label{eq:tdfc}
  \tilde\phi(t)=(1-R)\cdot\tilde\phi(t-T)+R\cdot\phi(t-T)  
\end{equation}
where $0<R\leq1$. Thus, whenever a dynamical system subject to
washout-filtered feedback control has a stable fixed point this fixed
point is also a (possibly unstable) fixed point of the
\emph{uncontrolled} dynamical system (the same is true for periodic
orbits of period $T$ using TDFC). This means that bifurcation analysis
could in principle be based on these non-invasive feedback control
techniques. One difficulty encountered, and intensively discussed in
the case of TDFC [\cite{NU98,JBORB97,HS05,FFGHS07}], is finding
reasonable conditions which guarantee that the non-invasive feedback
is actually able to stabilize a steady state of a dynamical
system. Typically, even if one has designed a feedback control loop
that is able to stabilize a steady state $\phi_*$ when one inserts
$\phi_*$ as the feedback reference signal (that is,
$\tilde\phi=\phi_*$) there is no guarantee that $\phi_*$ is also
stable if we replace the reference signal using the
recursion~\eqref{eq:tdfc} (which would correspond to TDFC).

\section{Pseudo-arclength continuation --- 
  vertically excited pendulum example}
\label{sec:parc}
\begin{figure}[t]
  \centering
  \includegraphics[width=0.9\columnwidth]{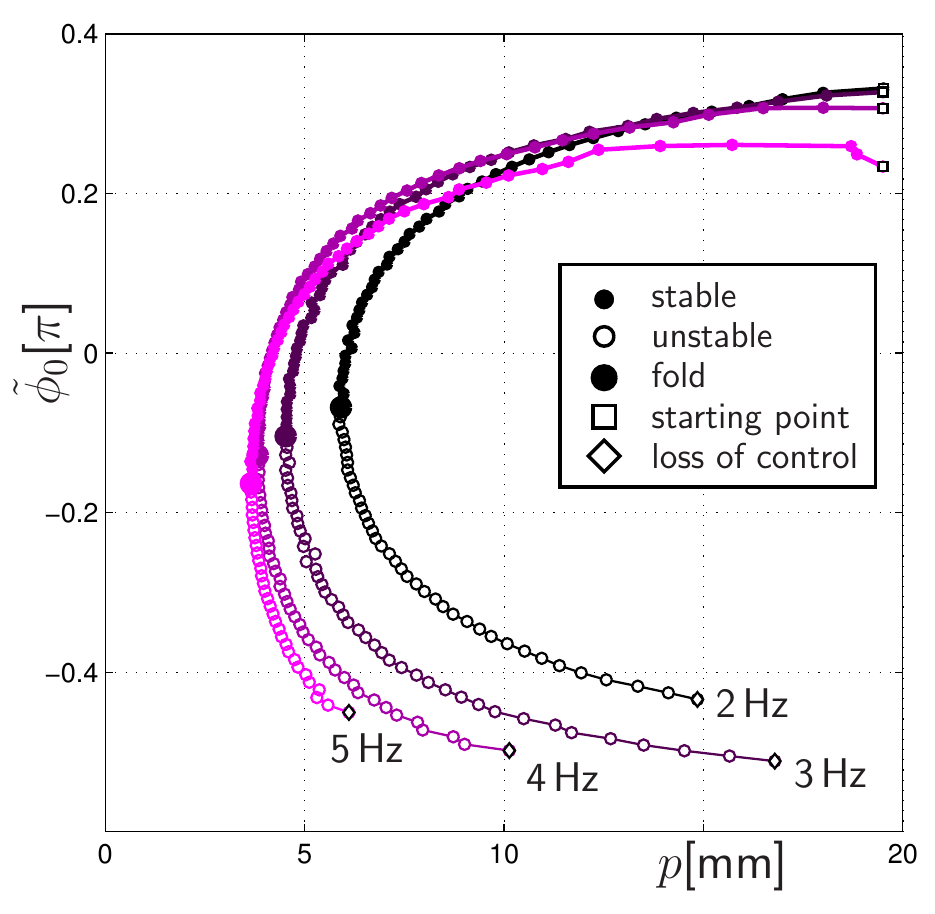}
  \caption{Bifurcation diagram for rotations of a parametrically
    excited pendulum at low forcing amplitude (experimental results
    from [\cite{SGNWK08}]). The $y$-axis shows the average phase
    $\tilde\phi_0$ as found by \eqref{eq:defpo}. Square marker: start
    of continuation; full circles: stable part of the branch of
    rotations; large circle: saddle-node bifurcation occurs in
    uncontrolled pendulum; hollow circles: unstable part of family of
    rotations; diamond: loss of control.}
  \label{fig:bif}
\end{figure}
If one shakes the pivot of a pendulum up and down harmonically with a
frequency $\omega$ and an amplitude $p$ then the pendulum can show
stable rotations for amplitudes larger than a certain critical
amplitude $p_0$. A simple model for the mechanical pendulum is
\begin{equation}
  \label{eq:pend}
  ml^2\ddot\theta+b\dot\theta+ml[g+\omega^2p\sin(\omega t)]
  \sin(\theta)=0\mbox{,}
\end{equation}
where $\theta$ is the angle, $m$ is the effective mass of the
pendulum, $l$ is its effective length, $b$ is a viscous damping
coefficient (typically small), $g$ is the acceleration due to gravity,
$\omega$ is the excitation frequency in rad$/$s, and $p$ is the
exciation amplitude in m. For the parameters used in Fig.~\ref{fig:bif}
the stable rotation at excitation amplitude $p\approx 2\,$cm (square
marker) is relatively easy to find by trial and error in a simulation
or an experiment because it has a large basin of attraction.

If one wants to find the minimal amplitude $p_0$ that supports
rotation one would in a simulation or experiment decrease the
excitation amplitude in small steps, always waiting until transients
decay after each parameter change. In this way one finds that the
stable periodic orbit ``disappears'' at (or, rather, slightly above)
the minimal amplitude $p_0$ corresponding to the saddle-node in
Fig.~\ref{fig:bif} as transients escape to another stable attractor,
for example, the hanging-down state. In a study of a model, such as
equation \eqref{eq:pend}, the alternative to this vary-and-wait
approach is a direct solution of \eqref{eq:pend} in
rotating coordinates $\phi=\theta-\omega t$:
\begin{equation}
  \label{eq:pendrot}
  ml^2\ddot\phi+b\dot\phi+b\omega+ml[g+\omega^2p\sin(\omega t)]
  \sin(\phi+\omega t)=0\mbox{,}
\end{equation}  with periodic
boundary conditions
\begin{equation}
  \label{eq:bc}
  \phi(T)-\phi(0)=0\mbox{,}\quad\dot\phi(0)-\dot\phi(T)=0\mbox{,}
\end{equation}
where $T=2\pi/\omega$ is the (known) period of the rotation. The
periodic boundary value problem \eqref{eq:pendrot},\,\eqref{eq:bc} is
nonlinear and is typically solved with a Newton iteration. The
advantage of this direct approach is that one can find periodic
rotations independent of their dynamical stability: the periodic
rotation undergoes a saddle-node bifurcation at the parameter value
$p_0$ such that the scenario looks as shown in Fig.~\ref{fig:bif}. The
Newton iterations for the nonlinear boundary value problem
\eqref{eq:pendrot},\,\eqref{eq:bc} finds both, dynamically stable and
unstable rotations in Fig.~\ref{fig:bif}.  Two difficulties for Newton
iterations are: (i) it converges only locally, that is, a good initial
guess is necessary; and (ii) the nonlinear problem is singular at the
saddle-node in Fig.~\ref{fig:bif}. Both problems can be overcome by
embedding the Newton iteration into a pseudo-arclength continuation
(see textbooks [\cite{D07,K04}]): one treats the bifurcation parameter
$p$ also as a variable, such that the solutions $(p,\phi(\cdot))$ of
\eqref{eq:pendrot},\,\eqref{eq:bc} form a curve ${\cal C}$ in the
space of all possible functions and parameters.  Furthermore, one
extends the nonlinear boundary value problem by the (scalar)
pseudo-arclength condition
\begin{multline}
  \label{eq:parc}
  \frac{1}{T}\int_0^T
  \dot\phi_\mathrm{tan}(t)[\dot\phi(t)-\dot\phi_\mathrm{old}]+
  \phi_\mathrm{tan}(t)[\phi(t)-\phi_\mathrm{old}(t)]\d t+\\
  +p_\mathrm{tan}[p-p_\mathrm{old}]=h\mbox{.}
\end{multline}
In \eqref{eq:parc} $(p_\mathrm{old},\phi_\mathrm{old}(\cdot))$ is the
previously found point on the curve ${\cal C}$,
$(p_\mathrm{tan},\phi_\mathrm{tan}(\cdot))$ is the unit tangent vector
to the curve ${\cal C}$ in this previous point and $h$ (a small
quantity) is the approximate distance between
$(p_\mathrm{old},\phi_\mathrm{old}(\cdot))$ and the desired solution
$(p,\phi(\cdot))$ of \eqref{eq:pendrot}--\eqref{eq:parc}.
Figure~\ref{fig:bif} shows a projection of this curve ${\cal C}$.
Using the pseudo-arclength extension \eqref{eq:parc} the nonlinear
boundary value problem is uniformly well-conditioned along the whole
curve ${\cal C}$ including the vicinity of the saddle-node at $p=p_0$.
Pseudo-arclength continuation is a useful (and, by now, well
established) tool in the numerical analysis of bifurcations because it
allows one to follow unstable parts of branches of periodic orbits and
also direct continuation of bifurcations (such as the saddle-node) in
more than one parameter. In this way one can construct maps in
the parameter space that help to classify for any given system its
possible equilibria and periodic orbits, and their bifurcations
[\cite{KOG07,K04}].

\section{Continuation in experiments}
\label{sec:expcont}
Extension of the pseudo-arclength continuation to experiments would
give experimenters the opportunity to study bifurcations in much
greater detail. For example, for the vertically excited pendulum with
rotations as shown in Fig.~\ref{fig:bif} a parameter study that simply
observes transients loses the stable periodic orbit already at
amplitudes $p$ significantly above $p_0$ (in the preparatory studies
for [\cite{SGNWK08}] at $p=0.8\,$cm for $\omega=3\,$Hz) due to
disturbances or insufficiently small parameter steps. Thus, it would
be difficult to establish that the loss of the stable periodic orbit
is indeed due to a saddle-node bifurcation. Moreover, the average
phase of the rotation with respect to the excitation changes
dramatically within a tiny parameter range: a stably rotating pendulum
points (nearly) upward whenever the pivot excitation reaches its
maximum (at $t=\pi/(2\omega)$) whereas close to the saddle-node
bifurcation the pendulum is nearly horizontal at time
$t=2\pi/\omega$. (Notice the extreme difference in the scaling of $x$-
and $y$-axis in Fig.~\ref{fig:bif} relative to measurement accuracy:
$\delta p\approx2\times10^{-4}$m, $\delta\phi\approx10^{-4}$rad.)  Due to
this sensitive dependence of the rotation on the parameter $p$, a
conventional experimental parameter study would also miss a
significant part of the upper stable part of the branch in
Fig.~\ref{fig:bif}.

\begin{figure}[t]
  \centering
  \includegraphics[width=\columnwidth]{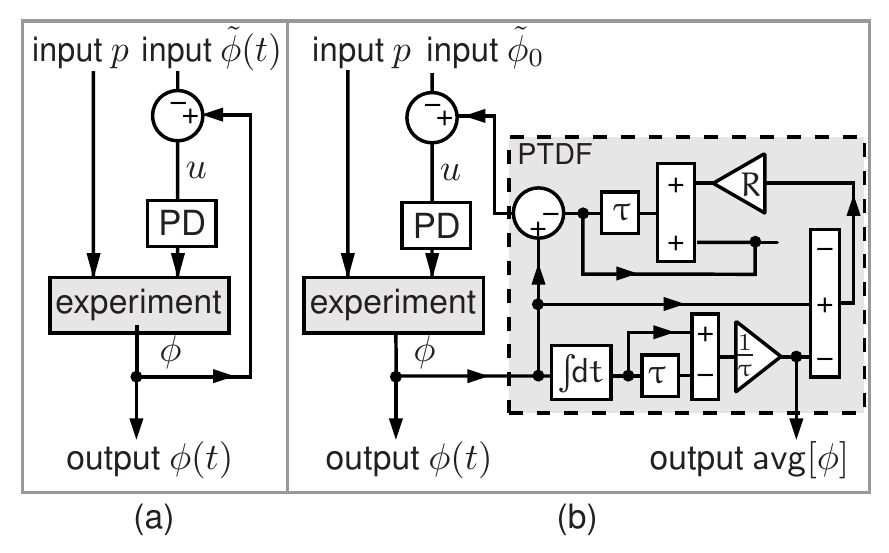}
  \caption{Pseudo block diagram for experiment with a tunable system
    parameter $p$ and a feedback loop: (a) standard feedback loop PD
    (proportional-plus-derivative) with periodic reference signal
    $\tilde \phi(t)$, (b) with projected time-delay (PTDF) block and
    constant scalar reference $\tilde\phi_0$. The square blocks with
    label $\tau$ delay the signal by $\tau=T=1\pi/\omega$.}
  \label{fig:flow}
\end{figure}
In [\cite{SK08b}] we proposed a method for continuing unstable
periodic orbits and bifurcations for experimental set-ups of the form
shown in Fig.~\ref{fig:flow}(a). The method assumes that a feedback
loop is implemented around the experiment (shown as a PD
(proportional-plus-derivative) controller in Fig.~\ref{fig:flow}(a))
which achieves stabilization in the following sense (formulated for
single-input-single-output feedback loops):
\begin{enumerate}
\item if the reference input $\tilde\phi(t)$ is identical to a
  periodic orbit $\phi_*(t;p_0)$ of the \emph{uncontrolled} experiment
  at parameter $p_0$ then $\phi_*(t;p_0)=\tilde\phi(t)$ is an
  exponentially \emph{stable} periodic orbit of the \emph{controlled}
  experiment. (It is a periodic orbit of the controlled experiment
  because the input $u$ into the PD control (see
  Fig.~\ref{fig:flow}(a)) vanishes for $\tilde\phi=\phi_*$.)
\item For parameter values $p\approx p_0$ and reference signals
  $\tilde\phi(t)\approx \phi_*(t;p_0)$ of period $T$ the output is
  asymptotically also periodic with period $T$ and the map
  \begin{equation}
    \label{eq:infmapspace}
    M_\infty:\R\times\cper([0,T];\R)\mapsto 
    \cper([0,T];\R)\mbox{,}
  \end{equation}
  defined by
  \begin{equation}
    \label{eq:eq:infmap}
    M_\infty(p;\tilde\phi(\cdot))(t)=\phi_\mathrm{asy}(t)
  \end{equation}
  is a locally well-defined and smooth map. The notation
  $\phi_\mathrm{asy}(t)$ refers to the (periodic) output of the
  controlled experiment with inputs $p$ and $\tilde\phi(t)$ after
  transients have died out, and $\cper([0,T];\R)$ is the space of all
  continuous real-valued periodic functions.
\end{enumerate}
We call the feedback loop stabilizing only if both conditions are
satisfied (possibly only locally). That a generic feedback loop can be
made stabilizing is a consequence of generic feedback stabilizability
of periodic orbits [\cite{NA92}]. The periodic orbits of the
uncontrolled experiment can be recovered as solutions of the nonlinear
fixed point problem
\begin{equation}
  \label{eq:fixedpointgen}
  M_\infty(p;\phi)=\phi\mbox{,}
\end{equation}
which, after discretization of $\phi$ (say into its $N$ first Fourier
coefficients), has one more variable ($N+1$) than equations ($N$) and
can be solved by a Newton iteration embedded into pseudo-arclength
continuation.  The Newton iteration defines a sequence of scalar
constant inputs $p$ and periodic inputs $\phi(\cdot)$, and the
residual required by the Newton iteration is the (periodic) asymptotic
limit of the control signal $u(\cdot)$ (see Fig.~\ref{fig:flow}(a)).
Since \eqref{eq:fixedpointgen} is an equation in the
infinite-dimensional space $\R\times\cper$ the accuracy of the result
(that is, how small $\|u\|_\infty$ can be made and, hence, how close
$\tilde\phi$ is to the unknown periodic orbit) depends on the choice
of $N$. (A larger $N$ means better approximation but also that a
larger number of repeated experiments with small input variations is
necessary.)

An alternative to the Newton iteration in the space $R\times\cper$ is
shown in Fig.~\ref{fig:flow}(b). The difference to the setup shown in
Fig.~\ref{fig:flow}(a) is the presence of an extra block, which we
called ``PTDF'' (for \emph{projected time-delayed feedback}) in
Fig.~\ref{fig:flow}(b), and which is inserted into the feedback loop
(inside the dashed rectangle). This block  implements the recursion
\begin{equation}
  \label{eq:rec}
  \tilde\phi(t)=(1-R)\tilde\phi(t-T)+R\cdot(\phi(t-T)-\avg[\phi(t-T-\cdot)])
\end{equation}
where $0<R\leq1$, and
\begin{equation}\label{eq:avg}
  \avg[\phi(t-\cdot)]=\frac{1}{T}\int_0^T\phi(t-s)\d s
\end{equation}
is the average of the signal $\phi$ over the past forcing period $T$.
The block feeds $\phi(t)-\tilde\phi(t)$ back into the feedback loop.
Thus, the input $u$ into the PD controller is
\begin{equation}\label{eq:uptdf}
  u(t)=\phi(t)-\tilde\phi(t)-\tilde\phi_0\mbox{.}
\end{equation}
This is a projected version of the extended time-delayed feedback
control (ETDFC; [\cite{P92,GSCS94}]): in the extreme case $R=1$ it
feeds back the difference between the output signal $\phi(t)$ and
$\phi(t-T)-\avg[\phi(t-T-\cdot)]+\tilde\phi_0$, which is the output
from one period ago but its average is shifted to the fixed input
$\tilde\phi_0$. For $R<1$ the function $\tilde\phi$ is a weighted sum
of the outputs from past periods (see \cite{GSCS94} for details).

The feedback controlled system in Fig.~\ref{fig:flow}(b) has the
property that, whenever the output $\phi$ of the experiment is
periodic (with the same period $T$ as the delay inside the block
``PTDF'' and the forcing), the output $\avg[\phi]$ is constant.
Furthermore, if the controlled system in Fig.~\ref{fig:flow}(b)
converges to a stable periodic motion $\phi_c(t)$ with period $T$ and
the limit of $\avg[\phi]$ is identical to its scalar input
$\tilde\phi_0$:
\begin{equation}\label{eq:avglim}
  \lim_{t\to\infty}\avg[\phi]=\avg[\phi_c]=\tilde\phi_0
\end{equation}
then $\phi_c(t)$ is a periodic orbit of the \emph{uncontrolled}
experiment. The control signal $u$ converges to zero if
\eqref{eq:avglim} is satisfied due to \eqref{eq:rec} and
\eqref{eq:uptdf}.

This implies that, if the projected ETDFC system shown in
Fig.~\ref{fig:flow}(b) has a two-parameter family of stable periodic
orbits in the $(p,\tilde\phi_0)$-plane of input parameters, we can
define the smooth map $M_1:\R^2\mapsto\R$ by
\begin{equation}
  \label{eq:M1def}
  M_1(p,\tilde\phi_0)=\lim_{t\to\infty}\avg[\phi(t-\cdot)]\mbox{,}
\end{equation}
and, whenever the input parameters satisfy the condition
\begin{equation}
  \label{eq:defpo}
  M_1(p;\tilde\phi_0)=\tilde\phi_0
\end{equation}
then the stable periodic orbit of the controlled system is identical to
a periodic orbit of the uncontrolled experiment. The curve shown in
Fig.~\ref{fig:bif} has been obtained as the curve of points in the
$(p;\tilde\phi_0)$-plane satisfying \eqref{eq:defpo}.

We note that \eqref{eq:defpo} is a scalar nonlinear equation in
contrast to the infinite-dimensional problem of the original algorithm
proposed in [\cite{SK08b}]. Here $M_1$ is obtained by setting the
input parameters in the controlled system in Fig.~\ref{fig:flow}(b),
waiting for transients to decay and then measuring the asymptotic
value of the scalar output $\avg[\phi]$, calling it
$M_1(p,\tilde\phi_0)$ and inserting it into \eqref{eq:defpo}.

Thus, the system shown in Fig.~\ref{fig:flow}(b) reduces the
infinite-dimensional nonlinear problem, as posed by the system in
Fig.~\ref{fig:flow}(a), to a scalar equation (which gives a
two-dimensional system after including the pseudo-arclength extension)
at the cost of the additional block in the feedback loop which has to
be evaluated in real-time in parallel to the experiment.

\section{Stability of the controlled system}
\label{sec:stab}

The remaining open questions are: suppose that the uncontrolled
experiment has a family of periodic orbits $\phi_*(t;p)$ and assume
that the feedback loop as shown in Fig.~\ref{fig:flow}(a) is
stabilizing for $\tilde\phi(t)=\phi_*(t;p)$.
\begin{enumerate}
\item\label{q:stab} When is the periodic orbit $\phi_*(t;p)$
of the corresponding projected time-delayed feedback controlled system
as shown in Fig.~\ref{fig:flow}(b) also stable for
$\tilde\phi_0=\avg[\phi_*]$?
\item\label{q:cond} What is the condition of the Jacobian of the
  reduced nonlinear problem? Sometimes, reducing the dimension of a
  nonlinear problem can cause a dramatic increase of its
  condition. (For example, if one reduces a periodic boundary-value
  problem to its corresponding fixed-point problem of the stroboscopic
  map.)
\end{enumerate}
Question \ref{q:stab} can be answered for a generalized version of the
projected time-delay block shown in Fig.~\ref{fig:flow}(b). Define the
(Fourier) spectral projections $P_N:C([0,T];\R)\mapsto\R^{2N+1}$ and $Q_N:\R^{2N+1}\mapsto C([0,T];\R)$:
\begin{displaymath}
  \begin{split}
    [P_Ny(\cdot)]_k&=\frac{1}{T}\int_0^Tb_k(2\pi s/T)y(s)\d s\mbox{,
      $k=-N\ldots N$}\\
    [Q_Nx](t)&=\sum_{k=-N}^Nx_kb_k(2\pi t/T)
  \end{split}
  \end{displaymath}
  where
  \begin{displaymath}
    \begin{split}
      b_k&=\sqrt{\smash[b]{2/T}}\cos(kt)\mbox{\ for $k<0$,}\\
      b_k&=\sqrt{\smash[b]{2/T}}\sin(kt)\mbox{\ for $k>0$, and}\\
      b_0&=\sqrt{\smash[b]{1/T}}\mbox{.}
    \end{split}
  \end{displaymath}
  The set-up in Fig.~\ref{fig:flow}(b)
  has to be generalized such that the input $\tilde\phi_0$ is not a
  constant but a periodic signal given by the combination of the
  harmonic oscillators of frequencies up to $N/T$ corresponding to the
  vector $x\in\R^{2N+1}$:
  \begin{displaymath}
    \tilde\phi_0(t)=[Q_Nx](t)\mbox{.}
  \end{displaymath}
  Furthermore, the
  projected time-delay block feeds 
  \begin{displaymath}
    \tilde\phi(t)=(1-R)\tilde\phi(t-T)+R[\phi(t-T)-P_N(\phi(t-T-\cdot))]
  \end{displaymath}
  back into the feedback loop and gives $P_N[\phi(t-\cdot)]$ as its
  $2N+1$-dimensional output (instead of the scalar $\avg[\phi]$).
  
  For this general projected time-delayed feedback control system the
  following holds: the periodic orbit $\phi_*(t;p)$ of the uncontrolled
  experiment is stable in the controlled system for sufficiently large
  $N$, sufficiently small $R$, and input parameters $\tilde\phi_0(t)$
  sufficiently close to $Q_NP_N\phi_*$. This generalized time-delayed
  feedback may be expensive (and impossible to perform in real-time)
  if the necessary $N$ is very large. It may also converge rather
  slowly if the necessary $R$ has to be chosen close to zero. 
  
  \begin{figure}[t]
    \centering
    \includegraphics[scale=0.48,angle=0]{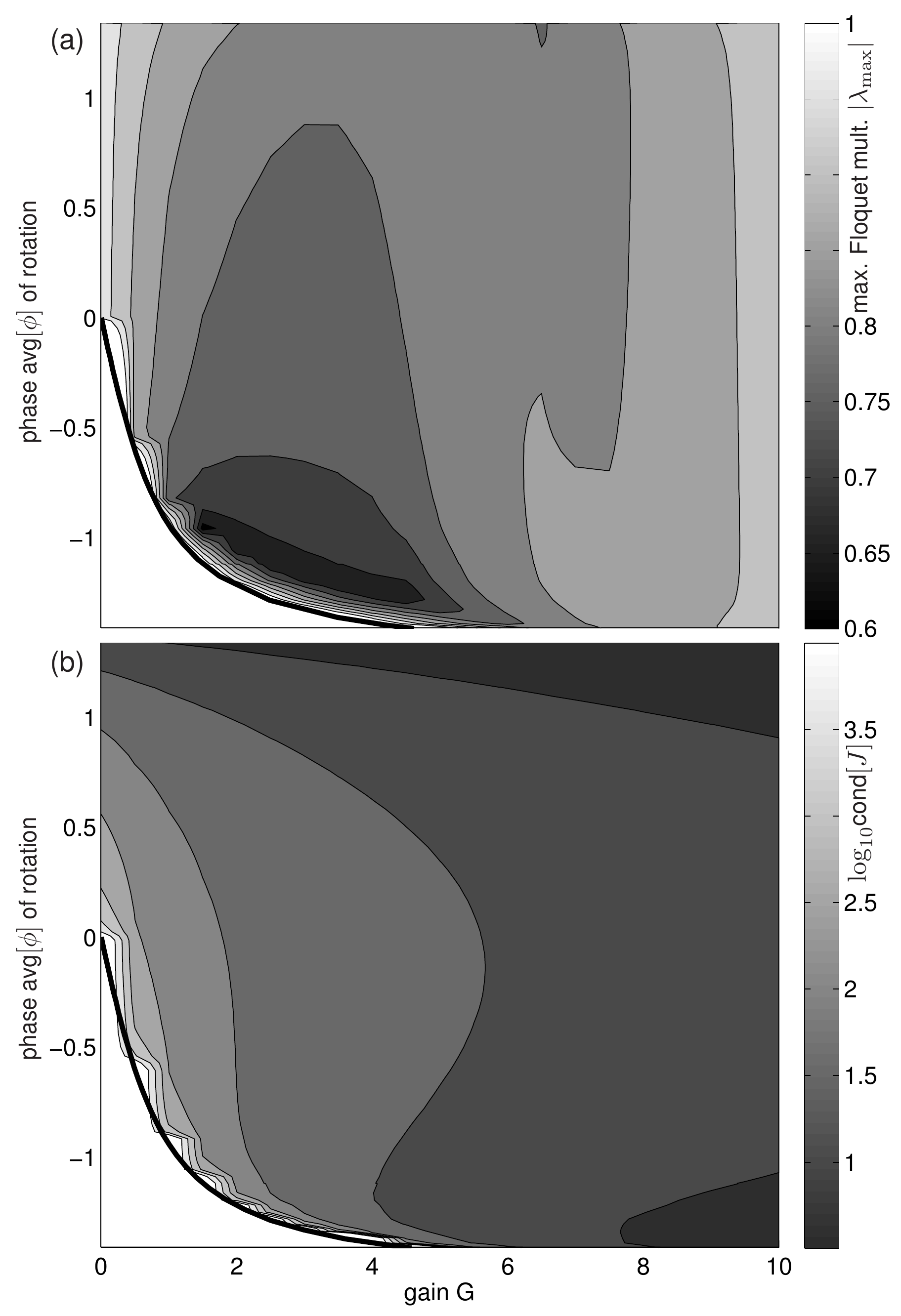}
    \caption{Stability (a) of the projected time-delayed feedback system
      shown in Fig.~\ref{fig:flow}(b) and condition of the resulting
      Jacobian matrix (b) for the vertically excited pendulum. The
      parameters in \eqref{eq:pendrot} are chosen to match the experimental
      data in [\cite{SGNWK08}; see caption of Fig.~\ref{fig:bif}].}
    \label{fig:stabcond}
  \end{figure}
  However, for a given experiment such as the rotations of the
  vertically excited pendulum studied in [\cite{SGNWK08}] the
  necessary $N$ may be very small and the permissible $R$ close to
  $1$. In fact, it turns out that for the study of rotations near the
  saddle-node bifurcation $N$ can be chosen equal to zero, which
  corresponds to the set-up in Fig.~\ref{fig:flow}(b). Moreover, the
  relaxation parameter $R$ can be chosen equal to $1$ which simplifies
  the input $u$ in the feedback loop to
  \begin{equation}\label{eq:utdf}
    u(t)=\phi(t)-\phi(t-T)-\avg[\phi(t-T-\cdot)]\mbox{.}
  \end{equation}
  Figure~\ref{fig:stabcond}(a) shows the stability of the periodic orbit
  of the uncontrolled system \eqref{eq:pendrot} when one applies the
  feedback loop shown in Fig.~\ref{fig:flow}(b) with $R=1$. The PD box
  in Fig.~\ref{fig:flow} has input $u$ and output
  \begin{equation}\label{eq:pd}
    PD[u]=- ml\,G\cdot[u+0.5\dot u]\mbox{,}
  \end{equation}
  which, in the simulation obtaining Fig.~\ref{fig:stabcond}, is added
  to the right-hand-side of the model \eqref{eq:pendrot} such that the
  overall controlled system in our simulations is
  \begin{equation}
    \label{eq:pdpend}  ml^2\ddot\phi+b\dot\phi+b\omega+ml[g+\omega^2p\sin(\omega t)]
  \sin(\phi+\omega t)=PD[u]
  \end{equation}
  where $PD[u]$ is defined by \eqref{eq:avg}, \eqref{eq:utdf} and
  \eqref{eq:pd}. This is an idealization of the experimental set-up in
  [\cite{SGNWK08}] where control had to be superimposed with the
  up-and-down excitation.  The $x$-axis in Fig.~\ref{fig:stabcond}
  shows the common factor $G$ of the control gains in the PD control
  \eqref{eq:pd}. The $y$-axis is the phase $\avg[\phi_*]$ of the
  periodic rotation of the uncontrolled system \eqref{eq:pendrot}.  As
  Fig.~\ref{fig:bif} shows, near the saddle-node bifurcation the
  family of periodic rotations of the uncontrolled pendulum
  \eqref{eq:pendrot} cannot be parametrized by the system parameter
  $p$ but by its phase $\avg[\phi_*]$. Each point in the plane in
  Fig.~\ref{fig:stabcond}(a) shows the largest Floquet multiplier of
  the periodic orbit $\phi_*(\cdot;p)$ as a periodic orbit of the
  controlled system \eqref{eq:pdpend}. The vertical line at $G=0$
  shows the stability of the periodic orbit without control: at $G=0$
  and $\avg[\phi_*]=0$ the dominant Floquet multiplier passes through
  $1$.  For varying $G$ this is a transcritical
  bifurcation. Increasing gains shift this loss of stability (black
  curve) along the originally unstable part of the branch toward lower
  $\avg[\phi_*]$.  Figure~\ref{fig:stabcond} shows that the controlled
  system is stable for sufficiently large gains $G$. The dominant
  Floquet multiplier does not decrease uniformly for increasing $G$
  because we chose the ratio between proportional and derivative term
  in the PD control \eqref{eq:pd} fixed (at $0.5$) and, thus,
  non-optimal. The stability chart looks similar for classical PD
  control using $u=\phi-\phi_*$ (that is, assuming that we knew the
  periodic rotations $\phi_*$ of the uncontrolled pendulum perfectly).

  The continuation procedure used in [\cite{SGNWK08}] to obtain the
  family of rotations in Fig.~\ref{fig:bif} starts from a stable
  rotation at $p\approx2\,$cm where one one can simply measure
  $\avg[\phi_*]$ and assign the initial input $\tilde\phi_0$ to this
  value. The initial unit tangent is
  $(p_\mathrm{tan},\tilde\phi_{0,\mathrm{tan}})=(-1,0)$. Then in each
  continuation step one performs a Newton iteration, running a
  sequence of controlled experiments as shown in
  Fig.~\ref{fig:flow}(b) for a sequence of inputs $(p,\tilde\phi_0)$
  as required by the Newton iteration and measuring the residual
  $r=(r_1,r_2)$ given by:
  \begin{equation}
    \label{eq:2dsys}
    \begin{split}
      r_1&=p_\mathrm{tan}\cdot(p-p_\mathrm{old})+
      \tilde\phi_{0,\mathrm{tan}}\cdot(\tilde\phi_0-\tilde\phi_{0,\mathrm{old}})-h\\
    r_2&=M_1(p;\tilde\phi_0)-\tilde\phi_0
    \end{split}
  \end{equation}
  where $(p_\mathrm{old},\tilde\phi_{0,\mathrm{old}})$ is the point
  previously found in the continuation, and
  $(p_\mathrm{tan},\tilde\phi_{0,\mathrm{tan}})$ is the (approximate)
  unit tangent to the solution curve in
  $(p_\mathrm{old},\tilde\phi_{0,\mathrm{old}})$. The Newton iteration
  is successful if the norm of the residual is smaller than a given
  tolerance ($5\times10^{-3}$ in [\cite{SGNWK08}]). The evaluation of
  $r_2$ requires running the controlled experiment until the
  transients have decayed. (How long this takes can be estimated from
  Fig.~\ref{fig:stabcond}(a).)

  This leads to the second question: how robustly does the Newton
  iteration converge? The convergence of the Newton iteration depends
  on the condition $\cond(J)=\|J\|\cdot\|J^{-1}\|$ of the Jacobian
  $J=[\partial r/\partial p,\partial r/\partial \tilde\phi_0]$, which
  is shown in Fig.~\ref{fig:stabcond}(b). We note that $\|J^{-1}\|$ is
  always $1$ for the system parameters in [\cite{SGNWK08}] (or very
  close to unity if the tangent is only approximate) because the rows
  of $J$ are orthogonal to each other by definition of the tangent
  $(p_\mathrm{tan},\tilde\phi_{0,\mathrm{tan}})$, and $\partial
  r_2/\partial p$ is large (due to the sensitive dependence of the
  rotation on the excitation amplitude and the moderate decay rate of
  the controlled system as shown in Fig.~\ref{fig:stabcond}(a)).
  Figure~\ref{fig:stabcond} provides evidence that the nonlinear
  system is uniformly well-conditioned for a wide range of gains. The
  only effect increasing the condition of the Jacobian is the loss of
  control when the gain becomes too small near the transcritical
  bifurcation (black curve in Fig.~\ref{fig:stabcond}(a) and (b)).

  \section{Conclusion}
  \label{sec:conc}

  We have analysed the stability and robustness of the experimental
  continuation of the periodic rotations of a pendulum through a
  saddle-node bifurcation performed in [\cite{SGNWK08}]. This analysis
  is important because the experiment relies on the asymptotic
  convergence of a feedback controlled experiment and Newton
  iterations, and in a real experiment one can never achieve that the
  input $u$ into the PD controller, defined by equation
  \eqref{eq:uptdf}, vanishes perfectly (due to disturbances and
  incomplete decay of transients, and a non-zero tolerance of the
  Newton iteration). We perform our analysis for a model equation.
  Because one is interested in the order of magnitude of condition
  numbers and decay rates. This approach is justified as long as the
  model is qualitatively correct (which is the case for a menchanical
  pendulum). We showed that both parts of the continuation process ---
  the experiment with its projected time-delayed feedback control and
  the Newton iteration in $\R^2$ for the inputs into the controlled
  experiment --- converge uniformly in the vicinity of the saddle-node
  bifurcation and on the unstable part of the family of rotations.

  \bibliography{delay} 
\end{document}